\documentclass[12pt]{article}
\usepackage[german,english]{babel}
\usepackage{amsmath,t1enc,times,theorem}
\usepackage{amsfonts}
\usepackage{multirow}

\newcommand{\1}{{\bf 1}}
\newcommand{\0}{{\bf 0}}
\newcommand{\R}{\mathbb{R}}

\newcommand{\N}{\mathbb{N}}
\newcommand{\E}{\mathbb{E}}
\newcommand{\prob}{\mathbb P}
\def\P{\prob}

\theoremstyle{plain} \newtheorem{Thm}{Theorem}
\theorembodyfont{\rmfamily} 
 
\newtheorem{Rem}{Remark}

\newenvironment{Prf}{\par \noindent {\bf Proof:}}{\\ $\Box$}

\begin{document}

\title{A semi-explicit density function for Kulkarni's bivariate phase-type distribution}
\author{Lothar Breuer}
\maketitle

\begin{abstract}
For a bivariate random variable $(Z_1,Z_2)$ having Kulkarni's bivariate phase-type distribution (see \cite{Ku89}), we derive a simple expression for the semi-explicit density function $\E (e^{-sZ_2} 1_{ \{ Z_1 \in dy \} })$. Some immediate consequences are presented as concluding remarks.
\end{abstract}

\section{Introduction}	\label{sec-ku}
Multivariate phase-type distributions have been a topic of research interest for quite some time. The first constructive proposal of such a class (henceforth denoted by MPH) can be found in \cite{AL84}. This class was later extended to MPH$^*$ in \cite{Ku89}, while the latest (and perhaps final) proposal of a definition (denoted by MVPH) is given in \cite{BN10}. All of the proposals carry their distinctive problems, be it that they seem too limited (as MPH) or that even elementary descriptions like distribution functions are not given explicitly (for MPH$^*$). 

The purpose of the present paper is the derivation of semi-explicit expressions for the density function of bivariate MPH$^*$ - distributed random variables. More exactly, for $(Z_1,Z_2) \in \text{MPH}^*$ we shall derive a simple expression for $\E (e^{-sZ_2} 1_{ \{ Z_1 \in dy \} })$, i.e.\ a density function for $Z_1$ joint with a Laplace transform for $Z_2$. As a univariate Laplace transform, this can be readily inverted to yield the bivariate density function $\P( Z_1 \in dy, Z_2 \in dx)$ for $x,y > 0$, see e.g.\ \cite{AW95}.

For ease of reference, we shall use the remainder of this introduction to restate the pertinent results in Kulkarni's construction of the class MPH$^*$. The main result along with some remarks are then presented in section 2.

Let ${\cal J} = ( J_t: t \geq 0)$ denote a Markov process on a finite state space $E' := \{ 1, \ldots, m+1 \}$ with $m \in \N$, having generator matrix
\[	\begin{pmatrix} 	Q & - Q \1 \\ \0 & 0
\end{pmatrix}
\]
where $Q$ is invertible, i.e.\ the states $i \in E := \{ 1, \ldots, m \}$ are transient. The initial distribution of ${\cal J}$ is denoted by $(\alpha, \alpha_{m+1})$ with $\alpha = (\alpha_1, \ldots, \alpha_m)$ and $\alpha_i := \P (J_0 = i)$ for $i \in E'$. We assume of course $\alpha_{m+1} < 1$. Let $R = (r_{ij})_{i \leq k, j \leq m}$ denote a reward matrix of dimension $k \times m$, with $r_{ij} \geq 0$ for all $i,j$. Write also $r_i(j) := r_{ij}$ whenever it is more convenient. Define the time of absorption of ${\cal J}$ by
\begin{equation}	\label{def-tau}
\tau := \min \{ t \geq 0: J_t = m+1 \}
\end{equation}
and further the random variables
\begin{equation}	\label{def-Z}
Z_i := \int_0^\tau r_i(J_t) dt
\end{equation}
for $i \in \{ 1, \ldots, k \}$. Then we say that $(Z_1, \ldots, Z_k) \in MPH^*$. The distribution of $(Z_1, \ldots, Z_k)$ shall be denoted by $MPH^*( \alpha, Q, R)$. To avoid trivial singularities later on, we assume that $\sum_{j=1}^m r_{ij} >0$, i.e.\ $\P(Z_i >0) >0$ for all $i \in \{ 1, \ldots, k \}$.

\section{The bivariate phase-type distribution}
From now on we set $k=2$, i.e.\ we consider bivariate distributions in MPH$^*$ only. The plan is the following: By theorem 1 in \cite{Ku89} the marginal distribution of $Z_1$ is phase-type. As indicated in remark 1 of \cite{Br12a}, it can be characterised in terms of the first passage times for a fluid flow. To be a bit more precise, let
\begin{equation} 	\label{def-fpt}
\tau(y) := \inf \{ t \geq 0: Y_t > y \}
\end{equation}
denote the first passage times for a suitable fluid flow model $({\cal J}, {\cal Y})$. Then 
\[	\P( Z_1 > y ) = \P_\alpha \left( \tau(y) < \tau | Y_0 = 0 \right) 
\]
where $\tau$ is the same as in (\ref{def-tau}) and $\P_\alpha$ denotes the conditional probability given that $\P( J_0 = i ) = \alpha_i$ for $i \in E$. Now we attach a phase-dependent time devaluation along the path of ${\cal Y}$ up to $\tau(y)$ to obtain an expression for 
\[	\E \left( e^{- s \int_0^{\tau(y)} r_2(J_s) ds} \right) 
\] 
which is the Laplace transform of $Z_2$ (with argument $s$) on the set of paths that satisfy $\tau(y) < \tau$, i.e.\ $Z_1 > y$. From here it is only a small step to obtain an expression for $\E ( e^{-s Z_2} 1_{ \{ Z_1 \in dy \} })$.

Two-dimensional fluid flow models have been analysed in detail in \cite{BR13}. We shall make use of some of the results therein, adapted to the question investigated here. In order to do so, we need to introduce some more notation. First we define the fluid flow models $({\cal J}, {\cal Y})$ and $({\cal J}, {\cal X})$ by 
\begin{equation}	\label{def-XY}
Y_t := \int_0^t r_1(J_s) \; ds \qquad \text{and} \qquad X_t := \int_0^t r_2(J_s) \; ds
\end{equation}
for all $t \geq 0$, where the phase process ${\cal J}$ is the same as in section \ref{sec-ku} and 
\[	r_1(m+1) := r_2(m+1) := 0
\]
Partition the set $E$ of transient states into $E=E_0 \cup E_+$, where
\[	E_0 := \{ i \in E: r_{1i} =0 \} \qquad \text{and} \qquad E_+ := \{ i \in E: r_{1i} > 0 \}
\]
According to this partition, write $Q$ and $\alpha$ in block form, i.e.\
\[	Q = \begin{pmatrix}	Q_{00} & Q_{0+} \\ Q_{+0} & Q_{++}	\end{pmatrix} \qquad \text{and} \qquad \alpha = (\alpha_0, \alpha_+)
\]
Further write $(\eta_0, \eta_+)' := \eta := -Q \1$. Finally, define the diagonal matrices 
\[	R_+ := diag( r_{1i}: i \in E_+),
\] 
\[ D_+ := diag( r_{2i}: i \in E_+) \qquad \text{and} \qquad D_0 := diag( r_{2i}: i \in E_0).
\]
Now we can state the main result:

\begin{Thm}
Let $(Z_1, Z_2) \sim MPH^*( \alpha, Q, R)$. Then 
\[	\E ( e^{-s Z_2} 1_{ \{ Z_1 =0 \} }) = \alpha_0 (s D_0 - Q_{00})^{-1} \eta_{0}
\]
for $s \geq 0$ and
\[	\E ( e^{-s Z_2} 1_{ \{ Z_1 \in dy \} }) = \alpha(s) e^{W(s) y} \eta(s) \; dy
\]
for $y > 0$ and $s \geq 0$, where
\begin{align*}
\alpha(s) &:= \alpha_0 (s D_0 - Q_{00})^{-1} Q_{0+} + \alpha_+ \\
W(s) &:= R_+^{-1} \left( (Q_{++} - s D_+) - Q_{+0} (Q_{00} - s D_0)^{-1} Q_{0+} \right) \\
\eta(s) &:= R_+^{-1}  \left( Q_{+0} (s D_0 - Q_{00})^{-1} \eta_0 + \eta_+ \right)
\end{align*}
\end{Thm}

\begin{Prf}
Due to the construction in (\ref{def-Z}) and (\ref{def-XY}), the representations $Z_1 = Y_\tau$ and $Z_2 = X_\tau$ hold, where $\tau$ is defined in (\ref{def-tau}). 

This means that on the set $\{ Z_1 =0 \}$, the phase process ${\cal J}$ lives only on $E_0$ before it gets absorbed. Define $\sigma := \min \{ t \geq 0: J_t \notin E_0 \}$. Clearly, $\sigma \leq \tau < \infty$ and $\{ \sigma = \tau \} = \{ Z_1 = 0 \}$. Thus
\[	\E ( e^{-s Z_2} 1_{ \{ Z_1 =0 \} }) = \E ( e^{-s \int_0^\tau r_2(J_s) \; ds} 1_{ \{ \sigma = \tau \} })
\]
Theorem 1 in \cite{BR13} states that
\[ \E ( e^{-s \int_0^t r_2(J_s) \; ds} 1_{ \{ t < \sigma < \tau \} }) = \alpha_0 e^{(Q_{00} - s D_0) t} \1
\]
for $s \geq 0$. Hence,
\[	\E ( e^{-s \int_0^t r_2(J_s) \; ds} 1_{ \{ \sigma = \tau \in dt \} }) = \alpha_0 e^{(Q_{00} - s D_0) t} \eta_0
\]
for all $t > 0$. Now integrating over $t \in ]0, \infty[$ yields the first statement. For the second statement, theorem 2 in \cite{BR13} states that
\[	\E \left( e^{-s \int_0^{\tau(y)} r_2(J_s) \; ds} 1_{ \{ \tau(y) < \tau \} } \right) = e^{W(s) y}
\]
for $s \geq 0$, where $\tau(y)$ is defined in (\ref{def-fpt}). Given our construction of ${\cal Y}$ and $Z_1$, this is equivalent to
\[	\E \left( e^{-s \int_0^{\tau(y)} r_2(J_s) \; ds} 1_{ \{ Z_1 > y \} } \right) = e^{W(s) y}
\] 
From here we obtain for small $h > 0$
\begin{multline*}
\E \left( e^{-s \int_0^{\tau(Z_1)} r_2(J_s) \; ds} 1_{ \{ y < Z_1 < y+h \} } \right) \\
= e^{W(s) y} \left( h R_+^{-1} \eta_+ + h R_+^{-1} Q_{+0} (s D_0 - Q_{00})^{-1} \eta_0 + o(h) \right)
\end{multline*}
and hence
\begin{align*}
\E ( \left. e^{-s Z_2} 1_{ \{ Z_1 \in dy \} } \right| J_0 = i) &=  e_i' e^{W(s) y} \eta(s) \; dy
\end{align*}
for $y > 0$ and ascending phases $i \in E_+$.  Considering all possible initial phases, we obtain by the same reasoning as for the first statement 
\[	\E ( e^{-s Z_2} 1_{ \{ Z_1 \in dy \} }) = \left( \alpha_0 (s D_0 - Q_{00})^{-1} Q_{0+} + \alpha_+ \right) e^{W(s) y} \eta(s) \; dy
\]
for $y > 0$, which is the second statement.  
\end{Prf}

\begin{Rem}
Let $Z=(Z_1, \ldots, Z_k) \in \text{MPH}^*$ with $k \geq 3$. According to theorem 6 in \cite{Ku89}, every pair $(Z_i, Z_j)$ with $i \neq j$ has a bivariate MPH$^*$ distribution. Thus we can use theorem 1 to determine the two-dimensional marginal distributions of a $k$-variate MPH$^*$ distribution.
\end{Rem}

\begin{Rem}
For $s=0$ we obtain the marginal distribution of $Z_1$, which is given as follows. Let $k := | E_0 |$ and $n := | E_+ |$, where $|M|$ denotes the cardinality of a set $M$. $Z_1$ has a PH($\beta, T$) distribution of order $n$ with
\[	\beta_i = \alpha_{k+i} + \alpha_0 (- Q_{00}^{-1}) Q_{0+} e_i
\]
for $i \in \{ 1, \ldots, n \}$ and
\[	\beta_{n+1} = \alpha_{m+1} + \alpha_0 (- Q_{00}^{-1}) \eta_0
\] 
The rate matrix $T$ is given by $T = W(0) = R_+^{-1} \left( Q_{++}  - Q_{+0} Q_{00}^{-1} Q_{0+} \right)$ such that
\begin{align*}
- T \1 &= - R_+^{-1} \left( Q_{++} \1 - Q_{+0} Q_{00}^{-1} Q_{0+} \1 \right)	\\
&= R_+^{-1} \left( \eta_+ + Q_{+0} \1 + Q_{+0} Q_{00}^{-1} (\eta_0 - Q_{00} \1) \right)	\\
&= R_+^{-1} \left( \eta_+ + Q_{+0} Q_{00}^{-1} \eta_0 \right)	\\
&= \eta(0)
\end{align*}
as to be expected.
\end{Rem}

\begin{Rem}
Theorem 4 in \cite{Ku89} states that the joint Laplace transform of $(Z_1, Z_2)$ is given by
\begin{equation}	\label{eq-ku}
\E ( e^{- s_1 Z_1} e^{- s_2 Z_2} ) = - \alpha ( \Delta - Q)^{-1} Q \1
\end{equation}
where $\Delta := diag( s_1 r_1(j) + s_2 r_2(j): j \in E)$. A relatively arduous way to arrive at this result is
\begin{align*}
\E ( e^{- s_1 Z_1} e^{- s_2 Z_2} ) &= \int_0^\infty e^{- s_1 y} \E ( e^{-s_2 Z_2} 1_{ \{ Z_1 \in dy \} }) \; dy + \E ( e^{-s_2 Z_2} 1_{ \{ Z_1 =0 \} }) \\
&= \left( \alpha_0 (Q_{00} - s_2 D_0)^{-1} Q_{0+} + \alpha_+ \right) \int_0^\infty e^{- s_1 y} e^{ W(s_2) y} \eta(s_2)  \; dy \\
& \qquad - \alpha_0 (Q_{00} - s_2 D_0)^{-1} \eta_{0} \\
&= - \left( \alpha_0 (Q_{00} - s_2 D_0)^{-1} Q_{0+} + \alpha_+ \right) (W(s_2) - s_1 I)^{-1} \eta(s_2)  \\
& \qquad - \alpha_0 (Q_{00} - s_2 D_0)^{-1} \eta_{0}
\end{align*}
First we observe that
\begin{align*}
W(s_2) - s_1 I &= R_+^{-1} \left( (Q_{++} - s_2 D_+) - Q_{+0} (Q_{00} - s_2 D_0)^{-1} Q_{0+} - s_1 R_+ \right) \\
&= R_+^{-1} \left( (Q_{++} - s_1 R_+ -s_2 D_+) - Q_{+0} (Q_{00} - s_1 R_0 - s_2 D_0)^{-1} Q_{0+} \right) 
\end{align*}
since $R_0 = \0$ by definition. To shorten notations, we write $W:= W(s_2) - s_1 I $. Further, we write
\[ R_+^{-1} \left( Q_{+0} (sD_0 - Q_{00})^{-1} \eta_0 + \eta_+ \right) = R_+^{-1} \left( - Q_{+0} (Q_{00} - s D_0)^{-1}, I \right) \begin{pmatrix} \eta_0 \\ \eta_+ \end{pmatrix}
\]
To arrive at (\ref{eq-ku}), we need to show that
\begin{align*}
( \Delta - Q)^{-1} &= \begin{pmatrix} (Q_{00} - s_2 D_0)^{-1} Q_{0+} \\ I \end{pmatrix} (-W)^{-1} R_+^{-1} \begin{pmatrix} - Q_{+0} (Q_{00} - s_2 D_0)^{-1} & I \end{pmatrix} \\
& \qquad + \begin{pmatrix} -(Q_{00} - s_2 D_0)^{-1} & \0 \\ \0 & \0 \end{pmatrix}
\end{align*}
In block form we can write
\[	( \Delta - Q) = \begin{pmatrix} s_2 D_0 - Q_{00} & Q_{0+} \\ Q_{+0} & s_1 R_+ + s_2 D_+ - Q_{++} \end{pmatrix}
\]
since $R_0 = \0$. Thus
\begin{align*}
( \Delta - Q) & \begin{pmatrix} (Q_{00} - s_2 D_0)^{-1} Q_{0+} \\ I \end{pmatrix} (-W)^{-1} \begin{pmatrix} - Q_{+0} (Q_{00} - s_2 D_0)^{-1} & I \end{pmatrix} \\
&= \begin{pmatrix} \0 \\ - R_+ W \end{pmatrix}  (-W)^{-1} R_+^{-1} \begin{pmatrix} - Q_{+0} (Q_{00} - s_2 D_0)^{-1} & I \end{pmatrix} \\
&= \begin{pmatrix} \0 & \0 \\ - Q_{+0} (Q_{00} - s_2 D_0)^{-1} & I \end{pmatrix} 
\end{align*}
and further
\[	( \Delta - Q) \begin{pmatrix} -(Q_{00} - s_2 D_0)^{-1} & \0 \\ \0 & \0 \end{pmatrix} = \begin{pmatrix} I & \0 \\  Q_{+0} (Q_{00} - s_2 D_0)^{-1} & \0  \end{pmatrix}
\]
Together this yields the desired result.
\end{Rem}

\begin{Rem} The most important ingredient to compute the covariance is $\E( Z_1 Z_2)$. Corollary 1 in \cite{Ku89} provides an iteration scheme to compute joint moments. An explicit formula is obtained via
\[ \E( Z_1 Z_2) = \left. \frac{d}{ds} \E( e^{-s Z_2} Z_1) \right|_{s=0}
\]
To this aim, 
\begin{align*}
\E( e^{-s Z_2} Z_1) &= \int_0^\infty y \E ( e^{-s Z_2} 1_{ \{ Z_1 \in dy \} }) = \alpha(s) \int_0^\infty y e^{W(s) y} \; dy \; \eta(s) 
\end{align*}
and
\begin{align*}
\int_0^\infty y e^{W(s) y} \; dy &= W(s)^{-1} \left[ y e^{W(s) y} \right]_{y=0}^\infty - W(s)^{-1} \int_0^\infty e^{W(s) y} \; dy = W(s)^{-2}
\end{align*}
yield
\begin{align*}
\E( Z_1 Z_2)  &= \left. - \frac{d}{ds} \alpha(s) W(s)^{-2} \eta(s)  \right|_{s=0} 
\end{align*}
This can be readily evaluated using the differentiation rule 
\[	\frac{d}{ds} (M(s)^{-1}) = M(s)^{-1} \frac{d}{ds} M(s) M(s)^{-1}
\]
for matrix-valued functions $M(s)$, see sections I.1.3-4 in \cite{Bo04}. 
\end{Rem}

\begin{Rem}
The special case MPH as described in \cite{AL84} is obtained as follows. Using the decomposition of state space $E$ and generator matrix $A$ as on p.692 therein, we can translate $E_+ = \Gamma_2^c$, $E_0 = \Gamma_2$, and
\[	Q_{++} = \begin{pmatrix} A^{(1,2)} & B^{(1)} \\ \0 & A^{(1)} \end{pmatrix}, \qquad Q_{+0} = \begin{pmatrix} B^{(2)} \\ \0 \end{pmatrix}, \qquad Q_{0+} = \0, \qquad Q_{00} = A^{(2)} 
\]
The construction in \cite{AL84} further specifies $R_+ = I$, $D_0 = I$, and 
\[	D_+ = \begin{pmatrix} I^{1,2} & \0 \\ \0 & \0 \end{pmatrix}
\] 
where $I^{1,2}$ denotes the identity matrix on $ \Gamma_1^c \cap \Gamma_2^c$. This yields $\alpha(s) = \alpha_+$,  
\[	W(s) = Q_{++} - s D_+ = \begin{pmatrix} A^{(1,2)} - s I & B^{(1)} \\ \0 & A^{(1)} \end{pmatrix}
\]
and
\[	\eta(s) = \begin{pmatrix} B^{(2)} (sI - A^{(2)})^{-1} (- A^{(2)} \1) \\ \0 \end{pmatrix} + \begin{pmatrix} -( A^{(1,2)} \1 + B^{(1)} \1 + B^{(2)} \1) \\ - A^{(1)} \1 \end{pmatrix} .
\]
\end{Rem}

\begin{Rem} If $r_{1j} > 0$ for all $j \in E$, then $E=E_+$, hence $\E ( e^{-s Z_2} 1_{ \{ Z_1 =0 \} }) =0$, and
\[	\E ( e^{-s Z_2} 1_{ \{ Z_1 \in dy \} }) = \alpha_+ e^{W(s) y} \eta_+
\]
where 
\[	W(s) = R_+^{-1} \left( Q_{++}  - s D_+ \right)
\]
for all $s \geq 0$. If further $r_{1i}=r_{2i}$ for all $i \in E$ with $r_{2i} > 0$ and $q_{ij} = 0$ for $r_{2i}=0$ and $r_{2j} > 0$, then we obtain the special case of the class MPH where $Z_1 \geq Z_2$ almost surely. This specifies to $\Gamma_2 = \0$, $\Gamma_1 = \{ i \in E: r_{2i} = 0 \}$, as well as
\[	A^{(1,2)} = \left( \frac{q_{ij}}{r_{1i}} \right)_{i,j \in \Gamma_1^c}, \quad B^{(1)} = \left( \frac{q_{ij}}{r_{1i}} \right)_{i \in \Gamma_1^c, j \in \Gamma_1} \quad \text{and} \quad A^{(1)} = \left( \frac{q_{ij}}{r_{1i}} \right)_{i,j \in \Gamma_1}.
\] 
\end{Rem}

\begin{Rem}
With no additional effort, the current framework can be extended to allow $r_{2i} < 0$ for some $i \in E$. One needs to take care of the range of $s$ for the Laplace transform $\E ( e^{-s Z_2} 1_{ \{ Z_1 \in dy \} })$ to converge (but there is such one, see lemma 2 in \cite{BR13}) or consider Fourier transforms. Then $Z_2$ has a so-called bilateral phase-type distribution, i.e.\ it is the mixture of two random variables $Z_2^+$ and $Z_2^-$ where $Z_2^+$ and $- Z_2^-$ have phase-type distributions. In particular, $Z_2$ may also assume negative values now. Theorem 2.3.2 in \cite{AR05} states that bilateral phase-type distributions are (weakly) dense in the class of all distributions on $\R$. For the marginal distribution of $Z_2$ see \cite{As04}, for more on bilateral phase-type distributions see \cite{AR05}.
\end{Rem}

\end{document}